\documentclass[preprint,showpacs,preprintnumbers,amsmath,amssymb]{revtex4}

\usepackage{graphicx}
\usepackage{dcolumn}
\usepackage{bm}

\begin{document}

\preprint{}

\title{Quasilocal angular momentum and center of mass in \\ general relativity}
\author{Po-Ning Chen}
 \affiliation{%
Columbia University, Department of Mathematics, 2990 Broadway, New York, New York 10027, USA\\
}%

\author{Mu-Tao Wang}
 \affiliation{%
Columbia University, Department of Mathematics, 2990 Broadway, New York, New York 10027, USA\\
}%
\author{Shing-Tung Yau}%
\affiliation{%
Harvard University, Department of Mathematics, One Oxford Street, Cambridge, Massachusetts 02138, USA\\
}%

\thanks{Part of this work was carried out while all three authors were visiting Department of Mathematics of National Taiwan University and Taida Institute for Mathematical Sciences in Taipei, Taiwan. P.-N. Chen is supported by NSF grant DMS-1308164, M.-T. Wang is supported by NSF grant DMS-1105483 and S.-T. Yau is supported by NSF grant  PHY-0714648.} 

\date{\today}

\begin{abstract}
For a spacelike 2-surface in spacetime, we propose a new definition of quasi-local angular momentum and quasi-local center of mass, as an element in the dual space of the Lie algebra of the Lorentz group. Together with previous defined quasi-local energy-momentum, this completes the definition of conserved quantities in general relativity at the quasi-local level. We justify this definition by showing the consistency with the theory of special relativity and expectations on an axially symmetric spacetime. The limits at spatial infinity provide new definitions for total conserved quantities of an isolated system, which do not depend on any asymptotically flat coordinate system or asymptotic Killing field. The new proposal is free of ambiguities found in existing definitions and presents the first definition that precisely describes the dynamics of the Einstein equation. 

\end{abstract}

\maketitle

\section{\label{sec:level1} Introduction}

By Noether's principle, any continuous symmetry of the action of a physical system corresponds to a conserved quantity. In special relativity, integrating the energy-momentum tensor of matter density against Killing fields in $\mathbb{R}^{3,1}$ gives the well-defined notion of energy momentum 4-vector, angular momentum, and center of mass. In attempt to generalize these concepts to general relativity, one encounters two major difficulties. Firstly,  gravitation does not have mass density. 
Secondly, there is no symmetry in a general spacetime. As such, most study of conserved quantities are restricted to isolated systems on which asymptotically flat coordinates exist at infinity. However, it has been conjectured \cite{Penrose} for decades that a quasi-local description of conserved quantities should exist, at least for energy-momentum and angular momentum. These are notions attached to a spacelike 2-surface $\Sigma$ in spacetime. In \cite{Wang-Yau1}, a new definition of quasi-local energy-momentum and quasi-local mass was proposed using isometric embeddings of the 2-surface into $\mathbb{R}^{3,1}$. The expression originated from the boundary term in the Hamilton-Jacobi analysis of gravitation action \cite{Brown-York2, hh}. To each pair of $(i, t_0^\nu)$ in which $i:\Sigma\hookrightarrow \mathbb{R}^{3,1}$ is an isometric embedding and $t_0^\nu$ is a future time-like unit Killing field, a canonical gauge (see (3) in \cite{Wang-Yau1}) is chosen and a
  quasi-local energy is assigned. The pair is considered to be a quasi-local observer and the quasi-local mass is obtained by minimizing the quasi-local energy seen among admissible $(i, t_0^\nu)$. A critical point of the quasi-local energy is an optimal isometric embedding \cite{Wang-Yau2}. In this letter, we use the optimal isometric embedding and the canonical gauge  to transplant Killing fields in $\mathbb{R}^{3,1}$ back to the surface of interest in the physical spacetime. In particular, this defines quasi-local angular momentum and quasi-local center of mass with respect to rotation Killing fields and boost Killing fields. We refer to \cite{Szabados} for earlier work on the definition of quasi-local angular momentum, notably \cite{Penrose2}. This new proposal is further applied to study asymptotically flat spacetime and to define new total conserved quantities on asymptotically flat initial data set. 
There are several existing definitions of total angular momentum and total center of mass such as the Arnowitt-Deser-Misner (ADM) angular momentum (\cite{Arnowitt-Deser-Misner,Ashtekar-Hansen,Regge-Teitelboim}) and the center of mass proposed by Huisken-Yau, Regge-Teitelboim, Beig-\'OMurchadha, Christodoulou, and Schoen \cite{Huisken-Yau, Regge-Teitelboim, Beig-Omurchadha, Christodoulou, Huang}. Unlike these definitions which rely on an asymptotically flat coordinate system or an asymptotic Killing field, the new definition is free of such ambiguities. In this Letter, we show that the new definition satisfies highly desirable properties, and fully captures the dynamics of the Einstein equation.

\section{\label{sec:level2} Definition and properties of quasi-local conserved quantities}
Suppose $\Sigma$ is a spacelike surface in a time orientable spacetime $N$, $u^\nu$ is a future-pointing timelike unit normal,
and $v^\nu$ is a spacelike unit normal with $u^\nu v_\nu=0$ along $\Sigma$. We recall the \textit{mean curvature vector field}
$h^\nu=-kv^\nu+pu^\nu $ in \cite{Wang-Yau1} and its companion $j^\nu=k u^\nu-pv^\nu $. Both are normal vector fields along $\Sigma$ which are defined independent of the choice of $u^\nu$ and $v^\nu$. We assume the mean curvature vector field is spacelike everywhere along $\Sigma$ and the norm of the mean curvature vector is denoted by $|H|=\sqrt{k^2-p^2}>0$. 
 $h^\nu$ and $j^\nu$ also define a connection one-form $\alpha_H$ of the normal bundle of $\Sigma$ by 
 \[\alpha_H=\frac{1}{k^2-p^2}\pi_{\alpha}^\beta(h^\nu\nabla_\beta j_\nu)\] where $\pi_\alpha^\beta=\delta_{\alpha}^\beta-u^\beta u_\alpha+v^\beta v_\alpha$ is the projection from the tangent bundle of $N$ onto the tangent bundle of $\Sigma$. We choose local coordinates  $\{u^a\}_{a=1,2}$ on $\Sigma$ and express this one-form as $(\alpha_H)_a$ and the induced Riemannian metric on $\Sigma$ as a symmetric $(0,2)$ tensor $\sigma_{ab}$. The definition of quasi-local conserved quantities depends only on the data $(\sigma_{ab}, |H|, (\alpha_H)_a)$, or the induced Riemannian metric and the mean curvature vector field. 

Consider a reference isometric embedding $i:\Sigma \hookrightarrow \mathbb{R}^{3,1}$ of $\Sigma$ so that the induced metric on the image surface is $\sigma_{ab}$. Suppose the mean curvature vector of the image surface is also spacelike  (this unnecessary assumption makes the exposition easier), we can similarly compute the mean curvature vector field on the image surface and obtain the correspond data $|H_0|$ and $(\alpha_{H_0})_a$.
 The definition of quasi-local energy, as in the Hamilton-Jacobi theory,  is with respect to a constant timelike
  unit vector $t_0^\nu$ in $\mathbb{R}^{3,1}$. Suppose the components of the isometric embedding $i$ is given by $(X^0, X^1, X^2, X^3)$, each as a smooth function on $\Sigma$. Let  $\eta_{\alpha\beta}$ be the Minkowski metric and $\tau=-t_0^\mu\eta_{\mu\nu} X^\nu$ be the time function on $\Sigma$ with respect to $t_0^\nu$. The quasi-local energy of $(\sigma_{ab}, |H|, (\alpha_H)_a)$ with respect to $(i, t_0^\nu)$ is given by 
  \begin{equation}\label{qle}\frac{1}{8\pi}\int_\Sigma \{( \cosh\theta_0|H_0|-\cosh\theta |H|)\sqrt{1+|\nabla\tau|^2}-[\nabla_a \theta_0-\nabla_a\theta+(\alpha_{H_0})_a-(\alpha_H)_a]\nabla^a\tau\}\end{equation}
 where $\theta_0=\sinh^{-1}\frac{-\Delta\tau}{|H_0|\sqrt{1+|\nabla\tau|^2}}$ and $\theta=\sinh^{-1}\frac{-\Delta\tau}{|H|\sqrt{1+|\nabla\tau|^2}}$. $\Delta\tau=\nabla^a\nabla_a \tau$ and $|\nabla\tau|^2=\sigma^{ab}\nabla_a\tau\nabla^a\tau$, where $\nabla_a$ is the covariant derivative of $\sigma_{ab}$. 
Expression \eqref{qle} is the same as (4) in \cite{Wang-Yau1} with $N_0=\sqrt{1+|\nabla\tau|^2}$ and $N_0^\nu=-\nabla^a\tau$. This  is derived from the boundary term of the surface Hamiltonian (see \cite{Brown-York2, hh}) of gravitation action in the canonical gauge condition (3) in \cite{Wang-Yau1}. 
  
It turns out the quasi-local energy is best described by $\rho $ and $j_a$ defined as follows: 
   \begin{equation} \label{rho} \begin{split}\rho &= \frac{\sqrt{|H_0|^2 +\frac{(\Delta \tau)^2}{1+ |\nabla \tau|^2}} - \sqrt{|H|^2 +\frac{(\Delta \tau)^2}{1+ |\nabla \tau|^2}} }{ \sqrt{1+ |\nabla \tau|^2}}. \end{split}\end{equation} and
\begin{equation} \label{j_a}
j_a=\rho {\nabla_a \tau }- \nabla_a [ \sinh^{-1} (\frac{\rho\Delta \tau }{|H_0||H|})]-(\alpha_{H_0})_a + (\alpha_{H})_a. \end{equation} 
Notice that $\theta_0-\theta=\sinh^{-1} (\frac{\rho\Delta \tau }{|H_0||H|})$ by the addition formula of the $\sinh$  function. 
In terms of these, the quasi-local energy is $\frac{1}{8\pi}\int_\Sigma (\rho+j_a\nabla^a\tau)$.
 
 We consider $(i, t_0^\mu)$ as a quasi-local observer and minimize quasi-local energy among all such observers. A critical point of the quasi-local energy satisfies the optimal isometric embedding equation :
 
\textit{Definition 1 \cite{Wang-Yau2, Chen-Wang-Yau1, Chen-Wang-Yau2} - An embedding $i:\Sigma\hookrightarrow \mathbb{R}^{3,1}$ satisfies the optimal isometric equation for $(\sigma_{ab}, |H|, (\alpha_H)_a)$ if the components of $i$, $X^0, X^1, X^2, X^3$, as functions on $\Sigma$, satisfy
 $\eta_{\mu\nu} \nabla_a X^\mu \nabla_b X^\nu=\sigma_{ab}$ and there exists a future unit timelike constant vector $t_0^\nu$ such that 
 $\tau=-t_0^\mu\eta_{\mu\nu} X^\nu$ satisfies
 \begin{equation}\label{optimal}\nabla^a j_a=0.\end{equation}}
 Such an optimal isometric embedding may not be unique, but it is shown in \cite{Chen-Wang-Yau2} that this is locally unique if $\rho>0$.
 
 The quasi-local mass and quasi-local energy-momentum 4-vector with respect to $(i, t_0^\nu)$ are $m= \frac{1}{8\pi}\int_\Sigma \rho $ and \[p^\nu= \frac{1}{8\pi}\int_\Sigma \rho \,t_0^\nu,\] respectively. 
 
Let $(x^0, x^1, x^2, x^3)$ denote the standard coordinate system on $\mathbb{R}^{3,1}$.  

\textit{Definition 2 -  Let $K_{\alpha\gamma}$ be an element of the Lie algebra of the Lorentz group with $K_{\alpha\gamma}=-K_{\gamma\alpha}$. Let  $K=K_{\alpha\gamma}\eta^{\gamma\beta} x^\alpha\frac{\partial}{\partial x^\beta}$ be the corresponding Killing field in $\mathbb{R}^{3,1}$. The conserved quantity corresponding to $(i, t_0^\nu, K)$ is $K_{\alpha\gamma}\Phi^{\alpha\gamma}$ where 
\begin{equation}\label{qlc_coordinate} \Phi^{\alpha\gamma}=-\frac{1}{8\pi} \int_\Sigma (\rho X^{[\alpha}  t_0^{\gamma]}+ j_a X^{[\alpha}  \nabla_a X^{\gamma]}). 
\end{equation}}
 
 For a spacelike 2-surface in $\mathbb{R}^{3,1}$, $\rho=0$ and $j_a=0$,  thus all quasi-local conserved quantities vanish with respect to its own isometric embedding.  By the definition, $p^\nu$ and $\Phi^{\alpha\gamma}$ transform equivariantly when the pair $(i, t_0^\mu)$ is acted by a Lorentz transformation. 
 When the optimal isometric embedding $i$ is shifted by a translation in $\mathbb{R}^{3,1}$ such that $X^\mu\mapsto X^\mu+b^\mu$ for some constant vector $b^\mu$, $\Phi^{\alpha\gamma}$ is changed by
 \[\Phi^{\alpha\gamma}\mapsto \Phi^{\alpha\gamma}-\frac{1}{2}b^\alpha p^\gamma+\frac{1}{2}b^\gamma p^\alpha.\]

For a surface of symmetry in an axially symmetric spacetime, the quasi-local angular momentum is the same as the Komar angular momentum, and the quasi-local center of mass lies on the axis of symmetry.
 
\section{\label{sec:level3} A conservation law}
Quasi-local conserved quantities satisfy a conservation law along timelike hypersurfaces described as follows.
The expression in \eqref{qle} can be written as the difference between a reference term and a physical term where the reference term is given by
\begin{equation}\label{qle.reference}\frac{1}{8\pi}\int_\Sigma \{ \cosh\theta_0|H_0|\sqrt{1+|\nabla\tau|^2}-[\nabla_a \theta_0+(\alpha_{H_0})_a]\nabla^a\tau\}.\end{equation}

On the image $i(\Sigma)$ of $i:\Sigma\hookrightarrow \mathbb{R}^{3,1}$, we choose the outward pointing spacelike unit normal $v_0^\nu$ such that $(t_0)_\nu v_0^\nu=0$. Let $u_0^\nu$ be the future pointing timelike unit normal of $i(\Sigma)$ such that $(u_0)_\nu v_0^\nu=0$. Extending $\Sigma$ alone the direction of $t_0^\nu$, we obtain a timelike hypersurface $\mathcal{C}$ of $\mathbb{R}^{3,1}$ whose spacelike unit outward normal is the extension of $v_0^\nu$. Let $\pi_{\mu\nu}$ be the conjugate momentum of $\mathcal{C}$. The expression \eqref{qle.reference} is the same as 
\begin{equation}\label{qle.reference2} \frac{1}{8\pi}\int_{i(\Sigma)} \pi_{\mu\nu} t_0^\mu u_0^\nu.\end{equation} 
Note that $i(\Sigma)$ is contained in $\mathcal{C}$ and $u_0^\nu$ is the normal of $i(\Sigma)$ in $\mathcal{C}$.  Since $\nabla^\mu\pi_{\mu\nu}=0$ and $t_0^\nu$ is Killing, we apply the divergence theorem to the portion of $\mathcal{C}$ that is bounded by $i(\Sigma)$ and a totally geodesic hyperplane orthogonal to $t_0^\nu$, and equate \eqref{qle.reference2} to the corresponding expression over the projection of $i(\Sigma)$ onto this hyperplane. The same procedure can be applied to the physical term by transplanting the Killing field $t_0^\nu$ back to the surface in the physical spacetime through the optimal isometric embedding and the canonical gauge. In a vacuum physical spacetime, $\nabla^\mu\pi_{\mu\nu}=0$ still holds, but the transplanted field may not be Killing and may not be tangent to the timelike hypersurface. However, these errors vanish at spatial infinity and we obtain a conservation law for the total mass and total angular momentum. This conservation law was first observed in  \cite{Brown-York2}. The novelty here is that this law is applied to both the physical and the reference spacetime.

\section{Conserved quantities at spatial infinity}

On an asymptotically flat initial data set $(M, g, k)$, we take the limit of quasi-local conserved quantities on coordinate spheres to define total conserved quantities. For each family of solutions $(i_r, t^\nu_0({r}))$ of optimal isometric embeddings for $\Sigma_r$ such that the isometric embedding $i_r$ converges to the standard embedding of a round sphere of radius $r$ in $\mathbb{R}^3$ when $r\rightarrow \infty$, we define:

 \textit{Definition 3 -Suppose $\lim_{r\rightarrow \infty} t^\nu_0({r})=t^\nu_0=L_0^{\,\,\, \nu}$ for a Lorentz transformation $L_\mu^{\,\,\,\nu}$.  Denote $L_{\mu\gamma}=L_\mu^{\,\,\,\nu}\eta_{\nu\gamma}$. The total center of mass of $(M, g, k)$ is defined to be 
\begin{equation}
C^i  =  \frac{1}{m}\lim_{r \to \infty}[ \Phi^{i\gamma}({r})L_{0\gamma}+\Phi^{0\gamma}({r})L_{i\gamma}], i=1, 2, 3
\end{equation} and the total angular momentum is defined to be
\begin{equation}
J_{i} =\lim_{r \to \infty} \epsilon_{ijk} [\Phi^{j\gamma}({r})L_{k\gamma}-\Phi^{k\gamma}({r})L_{j\gamma}], i,j,k=1, 2, 3 . 
\end{equation}}
$C^i$ corresponds to the conserved quantity of a boost Killing field, while $J_i$ corresponds to that of a rotation Killing field with respect to the direction of the energy-momentum 4-vector. 

\section{Finiteness of total conserved quantities}
In this section, we prove finiteness of the newly defined total conserved quantities for vacuum asymptotically flat initial data sets of order one.

\textit{Definition 4 - $(M, g, k)$ is \textit{asymptotically flat of order one} if there is a compact subset $C$ of $M$ such that 
$ M \backslash C $ is diffeomorphic to $\mathbb{R}^3 \backslash B$, and in terms of the coordinate system $\{x^i\}_{i=1, 2, 3}$ on $M\backslash C$, $g_{ij} = \delta_{ij}+ \frac{g_{ij}^{(-1)}}{r}+  \frac{g_{ij}^{(-2)}}{r^2}+ o(r^{-2})$ and $
k_{ij}  =   \frac{k_{ij}^{(-2)}}{r^2}  + \frac{k_{ij}^{(-3)}}{r^3}  + o(r^{-3})$,     
where $r=\sqrt{\sum_{i=1}^3 (x^i)^2}$.}

Transforming into spherical coordinates $(r, \theta, \phi)=(r, u^1, u^2)$, on each level set of $r$,  $\Sigma_r$, we can use $\{u^a\}_{a=1, 2}$ as coordinate system to express the geometric data we need in order to define quasi-local conserved quantities: 
\begin{equation}\label{expansion} 
\begin{split}
\sigma_{ab} & = r^2 \tilde \sigma_{ab}+ r \sigma_{ab}^{(1)} + \sigma_{ab}^{(0)}+ o(1) \\
|H| & = \frac{2}{r}+ \frac{h^{(-2)}}{r^2}+\frac{h^{(-3)}}{r^3} + o(r^{-3}) \\
\alpha_H & = \frac{\alpha_H^{(-1)} }{r}+ \frac{\alpha_H^{(-2)} }{r^2} + o(r^{-2}),
\end{split}
\end{equation} where $\tilde{\sigma}_{ab}$ is the standard round metric on $S^2$, and $\sigma_{ab}^{(1)}, \sigma_{ab}^{(0)}, h^{(-2)}, h^{(-3)}, \alpha_H^{(-1)}$, and $\alpha_H^{(-2)}$ are all considered as geometric data on $S^2$. 
It was proved in \cite{Chen-Wang-Yau1} that for such an  initial data set, there exists a unique family of optimal isometric embeddings $(i_r, t_0^\nu({r}))$ of $\Sigma_r$ such that the components of $i_r$ are given \begin{equation}\label{exp_opt}(0, r\tilde{X}^1, r\tilde{X}^2, r\tilde{X}^3)+o({r}),\end{equation} where $\tilde{X}^i, i=1, 2, 3$ denote the three coordinate functions on a unit 2-sphere in standard spherical coordinates. Similar expansions for $|H_0|$ and $\alpha_{H_0}$ can be obtained.

We recall from \cite{Wang-Yau3, Chen-Wang-Yau1} that the total energy-momentum vector $p^\nu$ satisfies  $\lim_{r \to \infty} t^\nu_0({r})=\frac{1}{m} p^\nu$ with $m=\sqrt{-p_\nu p^\nu}$, and the components are  given by
\begin{equation}\label{ADM}p^0= \frac{1}{8 \pi}\int_{S^2} (h_0^{(-2)}-h^{(-2)}) \text{ and } p^i= \frac{1}{8 \pi}\int_{S^2} \tilde{X}^i \widetilde{div}(\alpha_H^{(-1)}).\end{equation} It is shown to be the same as the ADM energy-momentum vector on the initial data set $(M, g, k)$ in \cite{Wang-Yau3}. 

Unlike translating Killing fields which define energy and linear momentum, the expression of boost and rotation Killing fields involves the coordinate functions. Therefore existing definitions of total angular momentum and total center of mass are in general infinite and  not well-defined unless additional conditions \cite{Ashtekar-Hansen, Regge-Teitelboim, Chrusciel2} are imposed at spatial infinity. 

Recall that an vacuum initial data set $(M, g, k)$ satisfies\begin{equation}\label{constraint}
R(g) + (tr_g k)^2-|k|_g^2 =0 \text{    and    } \nabla_g^i (k_{ij} - (tr_g k) g_{ij}) =0
\end{equation}
where $R(g)$ is the scalar curvature of $g_{ij}$.

 \textit{Theorem 1 -  The new total angular momentum and new center of mass are always finite on any vacuum asymptotically flat initial data set of order one.}

\textit{Proof.} By the expansions of geometric data \eqref{expansion} and the optimal isometric embedding \eqref{exp_opt}, the total center of mass and total angular momentum are finite if  
\begin{equation}\label{condition_converge}
\begin{split}
\int_{S^2} \tilde X^i (h_0^{(-2)}-h^{(-2)}) =0\text{ and } \int_{S^2}  \tilde X^i \left( \tilde{\epsilon}^{ab}\tilde{\nabla}_b(\alpha_H^{(-1)})_a \right)=0,
\end{split}
\end{equation} where $\tilde{\nabla}$ is the covariant derivative of $\tilde{\sigma}_{ab}$ and $\tilde{\epsilon}_{ab}$ is the area form on $S^2$.

Let $\Sigma_r$ be coordinate spheres of $(M, g)$ and $\hat h$ be the mean curvature of $\Sigma_r$ with respect to $g$. 
By the second variation formula of  area and the Gauss equation on $\Sigma_r$, we have
$ \partial_r \hat h = -f[\frac{R(g)}{2}-K+\frac{1}{2}(\hat h^2+A^2) ]- \Delta f, $
where $f$ is the length of $\frac{\partial}{\partial r}$, $K$ is the Gauss curvature, $A$ is the second fundamental form, and $\Delta$ is the Laplace operator of $\Sigma_r$. Expanding each term in the last equation in $r$, we solve for
$  \hat h^{(-2)} = -\tilde \sigma^{ab}\sigma^{(1)}_{ab} + \tilde \nabla^a \tilde \nabla ^b \sigma^{(1)}_{ab}  - \widetilde \Delta  (\tilde \sigma^{ab}\sigma^{(1)}_{ab})- \widetilde \Delta f^{(-1)} -2f^{(-1)}$, where $f=1+\frac{f^{(-1)}}{r}+o(r^{-1})$ and $\tilde{\Delta}$ is the Laplace operator of $S^2$. The decay condition on $k_{ij}$ implies 
$\hat h^{(-2)}=h^{(-2)}$ and thus $\int_{S^2} \tilde X^i h^{(-2)}=0$. $h_0^{(2)}$ can be solved from the optimal isometric embedding equation and $h_0^{(-2)}=-\frac{1}{2}\tilde{\sigma}^{ab}\sigma_{ab}^{(1)}-\tilde \nabla _e(\frac{1}{2} \tilde{\epsilon}^{ce} \tilde{\epsilon} ^{ad} \tilde\nabla _d\sigma^{(1)}_{ac} +\tilde{\epsilon}^{ce} F_c)$ for a one-form $F_c$ on $S^2$. Again, we obtain  $\int_{S^2} \tilde X^i h_0^{(-2)}=0$. On the other hand, by the vacuum momentum constraint equation, we derive $\tilde{\epsilon}^{ab}\tilde{\nabla}_b(\alpha_H^{(-1)})_a=  \tilde{\epsilon}^{ab}\tilde{\nabla}_b\tilde{\nabla}^c\pi^{(0)}_{ac}$ for a symmetric 2-tensor $\pi^{(0)}_{ac}$ on $S^2$, and thus the second equality in \eqref{condition_converge} also holds. This finishes the proof of Theorem 1.

The new total angular momentum vanishes on any hypersurface in $\mathbb{R}^{3,1}$: a property that is rather unique among known definitions. Indeed, it is shown in \cite{Chen-Huang-Wang-Yau} that there exists asymptotically flat hypersurface in $\mathbb{R}^{3,1}$ with finite non-zero ADM angular momentum.  In \cite{Chen-Wang-Yau3}, we show that the new total angular momentum is the same on any strictly spacelike hypersurface in the Kerr spacetime, a result of the conservation law discussed in \ref{sec:level3}.

\section{Conserved quantities and the vacuum Einstein equation}
Let $(M, g, k)$ be a vacuum initial data set. The vacuum Einstein equation is formulated as an evolution equation of the pair $(g(t), k(t))$ that satisfies $g(0)=g$, $k(0)=k$ and
\begin{equation} \label{evolution}
\begin{split}
\partial _t g_{ij} & =-2N k_{ij}+ \nabla_i \gamma_j+\nabla_j \gamma_i\\
\partial _t k_{ij} & = -\nabla_i\nabla_jN+N\left(R_{ij} + (tr_g k) k_{ij} - 2 k_{il} k^l \,_j\right)
\end{split}
\end{equation}
where $N$ is the lapse function and $\gamma$ is the shift vector.

\textit{Theorem 2- Suppose $(M, g, k)$ is an vacuum asymptotically flat initial data set of order one. Let $(M, g(t), k(t) )$ be the solution to the initial value problem $g(0)=g$ and $k(0)=k$ for the vacuum Einstein equation with lapse function $N=1+O(r^{-1})$ and shift vector $\gamma= \gamma^{(-1)}r^{-1}+O(r^{-2})$.
We have \[\begin{split}
\partial_t C^i (t)=  \frac{p^i}{p^0} \text{      and       }
\partial_t J_{i} (t) =   0
\end{split} \]
for $i=1, 2, 3$ where $p^\nu$ is the ADM energy momentum 4-vector of $(M, g, k)$.}

\textit{Proof.} By the expansions of geometric data \eqref{expansion} and the optimal isometric embedding \eqref{exp_opt}, it suffices to prove
\begin{equation}\label{center}\partial_t [\frac{1}{8\pi}\int_{S^2}\tilde{X}^i \rho^{(-3)}(t)]=\frac{p^i}{t_0^0} \text{    and    }\partial_t [\int_{S^2}(\tilde{X}^i\tilde{\nabla}_a\tilde{X}^j-\tilde{X}^i\tilde{\nabla}_a\tilde{X}^j) \tilde{\sigma}^{ab} j_b^{(-2)}(t)]=0.\end{equation}

In view of the definition of $\rho$ and the Einstein equation, $\int_{S^2} \tilde X^ i\partial_t \rho^{(-3)}(t)=\int_{S^2} \tilde X^ i\partial_t(\frac{h_0^{(-3)}(t)-h^{(-3)}(t)}{t_0^0})$. By applying the optimal isometric embedding and the vacuum constraint equation as in the proof of Theorem 1, together with the Einstein equation \eqref{evolution}, the first equality in \eqref{center} can be established.
As for the second equality, since $\tilde{\nabla}^a(\tilde{X}^i\tilde{\nabla}_a\tilde{X}^j-\tilde{X}^i\tilde{\nabla}_a\tilde{X}^j)=0$,
we can throw away two terms in $j_b$ whose leading terms are closed one-forms on $S^2$. The leading term of $\tau$ is $r\tau^{(1)}=-r \sum_k t_0^k \tilde{X}^k$ and
\[\partial_t[\int_{S^2}(\tilde{X}^i\tilde{\nabla}_a\tilde{X}^j-\tilde{X}^i\tilde{\nabla}_a\tilde{X}^j) \rho^{(-3)}(t)\tilde{\nabla}^a \tau^{(1)}]
=-t_0^j \partial_t[\int_{S^2} \rho^{(-3)}(t)\tilde{X}^i]+t_0^i \partial_t[\int_{S^2} \rho^{(-3)}(t) \tilde{X}^i)]\] where we use
$(\tilde{X}^i\tilde{\nabla}_a\tilde{X}^j-\tilde{X}^i\tilde{\nabla}_a\tilde{X}^j)\tilde{\nabla}^a \tilde{X}^k=\tilde{X}^i\delta^{jk}-\tilde{X}^j \delta^{ik}$. The last expression vanishes in view of the first equality in \eqref{center} and $p^\nu=mt_0^\nu$. The only term left is the time derivative of
$\int_{S^2}(\tilde{X}^i\tilde{\nabla}_a\tilde{X}^j-\tilde{X}^i\tilde{\nabla}_a\tilde{X}^j) \tilde{\sigma}^{ab} (\alpha_{H})_b^{(-2)}(t)$. This is shown to be independent of $t$ by the conservation law discussed in \ref{sec:level3}. This finishes the proof of Theorem 2. 

\noindent\textit{Remark - Assuming the total center of mass and angular momentum is finite initially, Theorem 2 holds under the weaker assumption $g= \delta + O(r^{-1})$ and $k= O(r^{-2})$, see \cite{Chen-Wang-Yau3}.}
\section{Properties of the new total conserved quantities}

1. The definition depends only on the geometric data $(g, k)$ and the foliation of surfaces at infinity, and in particular does not depend on an asymptotically flat coordinate system or the existence of an asymptotically Killing field. 

2. All total conserved quantities vanish on any spacelike hypersurface in the Minkowski spacetime, regardless of the asymptotic behavior.

3. The new total angular momentum and total center of mass are always finite on any vacuum asymptotically flat initial data set of order one.

4. The total angular momentum satisfies conservation law. In particular, the total angular momentum on any strictly spacelike hyperurface of the Kerr spacetime is the same. 

5. Under the vacuum Einstein evolution of initial data sets, the total angular momentum is conserved and the total center of mass obeys the dynamical formula $\partial_t C^i(t)=\frac{p^i}{p^0}$ where $p^\nu$ is the ADM energy-momentum four vector. 


\end{document}